\documentclass[12pt,leqno,fleqn]{amsart}
\usepackage{amsmath,amstext,amsthm,amssymb,amsxtra}
\usepackage[top=1.5in, bottom=1.5in, left=1.25in, right=1.25in]	{geometry}
\usepackage[normalem]{ulem}
\usepackage{txfonts} 
\usepackage[T1]{fontenc}
\usepackage{lmodern}
\usepackage{tikz}\usetikzlibrary{arrows}

\usepackage{euler}   

\makeatletter
\DeclareRobustCommand\widecheck[1]{{\mathpalette\@widecheck{#1}}}
\def\@widecheck#1#2{%
    \setbox\z@\hbox{\m@th$#1#2$}%
    \setbox\tw@\hbox{\m@th$#1%
       \widehat{%
          \vrule\@width\z@\@height\ht\z@
          \vrule\@height\z@\@width\wd\z@}$}%
    \dp\tw@-\ht\z@
    \@tempdima\ht\z@ \advance\@tempdima2\ht\tw@ \divide\@tempdima\thr@@
    \setbox\tw@\hbox{%
       \raise\@tempdima\hbox{\scalebox{1}[-1]{\lower\@tempdima\box
\tw@}}}%
    {\ooalign{\box\tw@ \cr \box\z@}}}
\makeatother

\usepackage{mathtools}
\mathtoolsset{showonlyrefs,showmanualtags}

\usepackage{hyperref} 
\hypersetup{
    colorlinks=true,       
    linkcolor=blue,          
    citecolor=magenta,        
    filecolor=magenta,      
    urlcolor=cyan           
}

\usepackage[msc-links]{amsrefs}

\theoremstyle{plain} 

\newtheorem{theorem}[equation]{Theorem}

\theoremstyle{definition}
\newtheorem{definition}[equation]{Definition}

\theoremstyle{remark}

\numberwithin{equation}{section}

\title[Positive Operators and Doubling Cubes]{Two Weight Inequalities  for Positive Operators: Doubling Cubes}

\author[W. Chen]{Wei Chen}
\address{ School of Mathematical Sciences, Yangzhou University, Yangzhou 225002, China; School of Mathematics, Georgia Institute of Technology, Atlanta GA 30332, USA}
\email {weichen@yzu.edu.cn}

\author[M. T. Lacey] {Michael T. Lacey}   

\address{ School of Mathematics, Georgia Institute of Technology, Atlanta GA 30332, USA}
\email {lacey@math.gatech.edu}
\thanks{The research of W. Chen is supported by the National Natural Science Foundation of China(11771379), the Natural Science Foundation of Jiangsu Province(BK20161326), and the Jiangsu Government Scholarship for Overseas Studies(JS-2017-228).
The research of M. T. Lacey is supported in part by grant from the US National Science Foundation, DMS-1600693 and the
Australian Research Council ARC DP160100153.}

\begin{document}
\begin{abstract}
For the maximal operator $ M $ on $ \mathbb R ^{d}$,  and $ 1< p , \rho < \infty $, 
there is a finite constant $ D = D _{p, \rho }$ so that this holds.  
For all  weights $ w, \sigma $ on $ \mathbb R ^{d}$,
 the operator  $ M (\sigma \cdot )$
is bounded from $ L ^{p} (\sigma ) \to L ^{p} (w)$ if and only the pair of weights $ (w, \sigma )$
satisfy the two weight $ A _{p}$ condition, and   this testing inequality holds: 
\begin{equation*}
\int _{Q}      M (\sigma \mathbf 1_{Q} ) ^{p}  \; d w \lesssim \sigma ( Q),
\end{equation*}
for all cubes $ Q$ for which there is a cube  $ P \supset Q$  satisfying $ \sigma (P) < D \sigma (Q)$, 
and $ \ell P = \rho \ell Q$. 
This  was recently proved by Kangwei Li and Eric Sawyer.
We give a  short proof, which is easily seen to hold for  several closely related operators.
\end{abstract}

	\maketitle
\tableofcontents

\section{Introduction} 

Our subject is the  two weight inequalities for the maximal function,
fractional integral transforms, and Poisson integrals.  For the purposes of this section, we will
focus on the maximal function.
A weight  $ w$ is a non-negative Borel measure on $ \mathbb R ^d$,
and given two weights $ w ,\sigma $  we say that  $ (w, \sigma ) \in A _{p}$ if the constant
\begin{equation} \label{e:Ap}
[w, \sigma ] _{p} = \sup _{Q}   \langle w  \rangle _Q ^{1/p} \langle \sigma  \rangle ^{1/p'}_Q , \qquad  p' = \tfrac p{p-1}.
\end{equation}
where here and throughout $  \langle w \rangle_Q = \lvert  Q\rvert ^{-1} \int _{Q} w \; dx $. With this notation the maximal function is
\begin{equation*}
M f = \sup _{ \textup{$ Q$} cube} \langle f \rangle_Q \mathbf 1_{Q}
\end{equation*}

The classical two weight inequality for the maximal function due to Sawyer \cite{MR676801} is below.
It shows that the inequality for the maximal function reduces to a testing inequality, for indicators of cubes.
\begin{theorem}\label{t:sawyer} For two weights $ (w, \sigma )$ we have the inequality
\begin{equation*}
\lVert M (\sigma f)\rVert _{L ^{p} (w)} \lesssim \lVert f\rVert _{L ^{p} (\sigma )}
\end{equation*}
if and only if the testing inequality below holds:
\begin{equation*}
\sup _{Q \;:\; \sigma (Q) >0} \sigma (Q) ^{-1/p} \lVert  \mathbf 1_{Q} M (\sigma \mathbf 1_{Q})\rVert _{L ^{p} (w)} < \infty .
\end{equation*}
\end{theorem}

Recent papers Li and Sawyer \cites{180904873,181111032} began a study of a weaker class of
testing inequalities in the two weight setting.  (Their papers include interesting motivation and background.)  They introduce four such conditions. The definition below is
  \emph{weaker than} their weakest condition: Test the maximal function on indicators of
cubes $ Q$ which have \emph{some} parent on which $ \sigma $ is doubling.

\begin{definition}\label{d:parent} Given two weights $ (w, \sigma )$, and $ 1< p , \rho , D< \infty $ we say that
$ (w, \sigma )$ satisfy a $ (p , \rho , D)$ parent doubling testing condition if there is a positive finite constant $ \mathfrak P = \mathfrak P _{\rho , D} = \mathfrak P (w, \sigma ,d,  p, \rho ,D)$ so that we have
\begin{equation}\label{e:parent1}
\lVert  \mathbf 1_{Q} M (\sigma \mathbf 1_{Q})\rVert _{L ^{p} (w)}\leq \mathfrak P \sigma (Q) ^{1/p},
\end{equation}
for every cube $ Q$ for which there is a second cube $ P \supset    Q$, with $   \ell P  \geq    \rho   \ell Q$,
and $ \sigma (P) \leq D \sigma (Q)$.
\end{definition}

Above $ \ell Q = \lvert  Q\rvert ^{1/d} $ is the side length of $ Q$.
The case of $ \rho =p=2$  below is (just a little stronger than)
the main result of Li and Sawyer \cite{181111032}.

\begin{theorem}\label{t:parent} Let   $ 1< p , \rho <\infty$.  There is a constant $ D = D _{d, p , \rho }$ so that
for any pair of weights $ (\sigma ,w)$ we have
\begin{equation*}
\lVert M (\sigma \cdot )\rVert _{L ^{p} (\sigma ) \to L ^{p} (w)}
\simeq    [w, \sigma ]_p + \mathfrak P _{\rho , D}.
\end{equation*}
\end{theorem}

Our proof is short and   readily adapts to
several closely related operators, as indicated in the concluding section.

\section{Proof} 
The norm bound on $ M$ easily implies the two weight $ A_p$ condition on the weights, as well as the testing condition for all cubes, not just those with a doubling parent.
The content of the Theorem is that the reverse implication holds.

Our theorem only claims that there is a sufficiently large doubling parameter $ D$ which
can be used for all pairs of weights $ (w, \sigma )$.
 Below, we will consider values of $ 1< \rho \leq 2 $.
 For integers $ n = 3, 4 ,\dotsc, $, and choices of $ n-1< \rho \leq n$, the argument proceeds by replacing
 the dyadic grids introduced below by $ n$-ary grids. We omit the details.

By a \emph{dyadic grid} we need a collection $ \mathcal D$ of cubes in $ \mathbb R ^d$
for which  (a) if $ P, Q \in \mathcal D$, then $ P\cap Q $ is empty, $ P$ or $ Q$,
and (b) for all integers $ k$, the cubes $ \{Q \in \mathcal D \;:\; \lvert  Q\rvert = 2 ^{dk} \}$
partition $ \mathbb R ^d$. Associated to the grid $ \mathcal D$ is the maximal function
\begin{equation*}
M _{\mathcal D} f = \sup _{Q\in \mathcal D} \langle f \rangle_Q \mathbf 1_{Q} .
\end{equation*}

As is well known, there are a finite number of grids $ \mathcal D_1 ,\dotsc, \mathcal D _{3 ^{d}}$
for which
\begin{equation} \label{e:MD}
M f \lesssim \sup _{1\leq j \leq 3 ^{d}}  M _{\mathcal D_j} f  .
\end{equation}
Set $ D =2 ^{d \frac{p+1} {p-1}}$. It suffices to  show that under the two weight $ A_p$ and  $ (p, 2 , D)$
parent testing condition, for any dyadic grid $ \mathcal D$, the maximal function $ M _{\mathcal D} (\sigma \cdot )$ is bounded from $ L ^{p} (\sigma )$ to $ L ^{p} (w)$.

Sawyer's Theorem  \ref{t:sawyer} holds for $ M _{\mathcal D}$.
Namely, it suffices to show that for any cube $ Q_0\in \mathcal D$, we have
\begin{equation} \label{e:PT}
\int _{Q_0}  M _{\mathcal D} ( \sigma \mathbf 1_{Q_0} ) ^{p} \; dw
 \lesssim ([w, \sigma ]_p ^{p} + \mathfrak P  ^{p}) \sigma (Q_0).
\end{equation}
We are free to restrict the supremum  to the collection of cubes
$ \mathcal Q = \{ {Q\in \mathcal D \;:\;  Q \subset Q_0}\}$ of cubes contained in $ Q_0$.

Partition $ \mathcal Q  $ into four subcollections using these definitions.
\begin{itemize}\setlength\itemsep{1em}
\item  (Testing Collection)   Let $ \mathcal T ^{\ast} $ be the maximal elements   $ Q\in \mathcal D$
with $ Q\subset Q_0$  so that the testing inequality \eqref{e:parent1} holds.
Set $ \mathcal T_Q = \{P\in \mathcal Q \;:\; P\subset Q\}$, for $ Q\in \mathcal T ^{\ast} $.  And set $ \mathcal T = \bigcup _{Q\in \mathcal T ^{\ast} } \mathcal T_Q$.

\item   (The Top)  Let $ \mathcal U = \{Q\in \mathcal Q \setminus \mathcal T \;:\;
2 ^{k}\ell Q \geq \ell Q_0\}$.
We choose $ k$ large enough that $ 2 ^{dk} k ^{-p} >1$.
These are the cubes which are close to the top cube  $ Q_0$.

\item  (Small $ A_p$ Cubes) Let $ \mathcal A $ be those cubes $ Q\in \mathcal Q \setminus (\mathcal T \cup \mathcal U ) $ such that
\begin{equation}\label{e:small}
\langle \sigma  \rangle_Q ^{1/p'} \langle w \rangle_Q ^{1/p} \leq \frac{[w, \sigma ] _{p} } {  \log  \ell   Q_0/  \ell   Q    } .
\end{equation}
That is, the local $ A_p$ constant at $ Q$ is very small.

\item (Remaining Cubes)  Let $ \mathcal R = \mathcal Q \setminus (\mathcal T \cup \mathcal U \cup \mathcal A)$.
\end{itemize}

We show that  the maximal function over each collection satisfies the testing inequality \eqref{e:PT}.
The Testing Collection is very easy:
\begin{align*}
\int _{Q_0}  \sup _{Q\in \mathcal T}  \langle \sigma  \rangle_Q^{p} \mathbf 1_{Q} \; dw
& \leq \sum_{Q \in \mathcal T ^{\ast} }  \int _{Q}  M (\sigma \mathbf 1_{Q}) ^{p}\; dw
\\
& \leq \mathfrak P _{ 2 ,D} ^{p} \sum_{Q \in \mathcal T ^{\ast} } \sigma (Q)
\leq  \mathfrak P  ^{p} \sigma (Q_0).
\end{align*}

The Top Collection $ \mathcal U$   has at most  $ 2 ^{1+d(k+1)}$ elements, and we just use the $ A_p$ condition
to see that
\begin{align*}
\int _{Q_0}  \sup _{Q\in \mathcal U}  \langle \sigma  \rangle_Q^{p} \mathbf 1_{Q} \; dw
& \leq
\sum _{Q\in \mathcal U}   \langle \sigma  \rangle_Q  ^{p} \langle w \rangle _{Q}
\\
& \leq  [w, \sigma ] _{p} ^{p} \sum _{Q\in \mathcal U}   \sigma (Q) \lesssim [w, \sigma ] _p ^{p} \sigma (Q_0).
\end{align*}
The implied constant depends upon $k$, but that is a fixed integer.

The Small $ A_p$ Cubes are also trivially sum up, using the condition in \eqref{e:small}.
\begin{align*}
\int _{Q_0}  \sup _{Q\in \mathcal A}  \langle \sigma  \rangle_Q^{p} \mathbf 1_{Q} \; dw
& \leq
\sum _{Q\in \mathcal A}   \langle \sigma  \rangle_Q  ^{p} \langle w \rangle _{Q}
\\
& \leq  [w, \sigma ] _{p} ^{p} \sum _{Q\in \mathcal A} \frac{\sigma (Q)}
 { [ \log   \ell   Q_0/  \ell   Q ] ^{p}  } \lesssim [w, \sigma ] _p ^{p} \sigma (Q_0).
\end{align*}

\smallskip
Thus, the core of the argument is control of the Remaining Cubes, $ \mathcal R$. Indeed, we claim that this
collection is empty, since a cube that has a large local $ A_p$ product is also approximately doubling.

Suppose $ \mathcal R \neq \emptyset $. Thus, there is a cube $ Q \subset Q_0$,
which satisfies $ \ell Q < 2 ^{-k} \ell Q_0 $,   fails \eqref{e:small}, and \emph{no ancestor of $ Q$ also  contained inside of $ Q_0$, has a doubling parent.}  The last condition is  very strong.

Let $ Q ^{(1)}$ be the $ \mathcal D$-parent of $ Q$, and let $ Q ^{(k+1)} = (Q ^{(k)}) ^{(1)}$.
Define integer $ m$ by  $ Q_0 = Q ^{(m)}$.
For any integer $ 0\leq k < m$, we
necessarily have $   \sigma (Q ^{(k+1)})\geq D  \sigma (Q^{(k)})$, since $ Q ^{(k+1)}$ is a $ \rho $-parent of
$ Q ^{(k)}$.
That is, $ \sigma (Q_0) \geq D ^{m} \sigma (Q)$.
From this, we see that $ m$ cannot be very large.
\begin{align*}
[w, \sigma ]_p ^{p}  & \geq \langle \sigma  \rangle_ {Q_0} ^{p-1} \langle w \rangle _{Q_0}
\\
& \geq   D ^{m (p-1)} \Bigl[\frac{\sigma (Q)} { \lvert  Q ^{(m)}\rvert } \Bigr] ^{p-1} \frac{w (Q)} {\lvert  Q ^{(m)}\rvert }
\\
& \geq   [D/ 2 ^{dp'}] ^{m (p-1)}   \langle \sigma  \rangle_Q ^{p-1} \langle w \rangle_Q
\\
& \geq [w, \sigma ]_p ^{p}  [D/ 2 ^{dp'}] ^{m (p-1)}  m ^{- p} \geq [w, \sigma ]_p ^{p} 2 ^{dm} m ^{-p}.
\end{align*}
The constants are explicit, and the last line follows by choice of $ D$.
We see that $ m < k$.  That is, the cube is in the collection $ \mathcal U$, which is a contradiction.

\section{Complements} 

\newcounter{para}
\newcommand\mypara{\par\smallskip\refstepcounter{para}\textbf{\thepara.}\space}

\mypara The conditions in Theorem~\ref{t:parent} can be strengthened by adding the
condition that the doubling cubes satisfy $ \sigma (\partial D) =0$.
This is accomplished by selection of \emph{random grids}. The discussion needed is given by
Li and Sawyer \cite{181111032}*{\S2}, and we omit the details.  Similar comments apply to the
extensions we mention below.

\mypara Theorem~\ref{t:parent} has a straight forward extension to fractional maximal functions.

\mypara The method of proof easily extends to other operators which are well approximated by
dyadic grids.  One of these is the Poisson integral given by
\begin{equation*}
P f (x,t) = \int \frac{t} {  (t ^2 + \lvert  x-y\rvert ^2  ) ^{\frac{d+1}2} } f (y) \; dy , \qquad t>0
\end{equation*}
Given weights $ \sigma $ on $ \mathbb R ^{d}$ and $ w$ on the upper half space $ \mathbb R ^{d+1} _{+}$,
we remark that the role of cubes in $ \mathbb R ^{d+1} _{+}  $ are played by \emph{Carleson cubes},
namely $ \tilde Q =  Q \times [0, \ell Q)$, for $ Q\subset \mathbb R ^{d}$.  The definition of the two weight $ A_p$ condition is then
\begin{equation}  \label{e:App}
[w, \sigma ] _{p} = \sup _{Q}  \langle w \rangle _{\tilde Q} ^{1/p} \langle \sigma  \rangle_Q ^{1/p'}.
\end{equation}

Using similar methods, one can prove this version of the Poisson two weight theorem of Sawyer \cite{MR719674}.  We single out this statement since the Sawyer's Poisson theorem is an important ingredient of the two weight inequality for the Hilbert transform \cites{MR3285857,MR3285858}.

\begin{theorem}\label{t:poisson} Let $ 1< p , \rho < \infty $.  There is a $ D>1$ so that this holds.
Let $ w$ be a weight on $ \mathbb R ^{d}_+ $, $ \sigma $ on $ \mathbb R ^{d}$.  These conditions are necessary and sufficient for $ P (\sigma   \cdot) $ to map $ L ^{p} (\mathbb R ^{d}, \sigma   )$
to $ L ^{p} (\mathbb R ^{d+1}_+, w ) $.   There is a finite constant $ \mathfrak P$ so that
\begin{enumerate}
\item   The $ A_p$ condition \eqref{e:App} holds.

\item  If $ Q\subset \mathbb R ^{d}$ is a cube for which $ \sigma (\rho Q) < D \sigma (Q)$,
then,   $ \lVert \mathbf 1_{\tilde Q} P (\sigma \mathbf 1_{Q})\rVert _{L ^{p} (\mathbb R ^{d+1}_+, w)} \leq \mathfrak P \sigma (Q) ^{1/p}$.

\item If $ Q\subset \mathbb R ^{d}$ is a cube for which $ w (  \widetilde {\rho Q}) < D \sigma (\tilde Q)$,
then,   $ \lVert \mathbf 1_{  Q} P ^{\ast} (w \mathbf 1_{\tilde Q})\rVert _{L ^{p'} (\mathbb R ^{d}, \sigma )} \leq \mathfrak P w (\tilde Q) ^{1/p'}$.

\end{enumerate}

\end{theorem}

Let us briefly indicate the proof. For any dyadic grid $ \mathcal D$, we can define
\begin{equation*}
P _{\mathcal D} f = \sum_{Q\in \mathcal D} \langle f \rangle _{Q} \mathbf 1_{\tilde Q} .
\end{equation*}
One has $ P _{\mathcal D} f \lesssim P f $, but also an analog of \eqref{e:MD} holds.
That is, there are finitely many dyadic grids
 $ \mathcal D_1 ,\dotsc, \mathcal D _{3 ^{d}}$
for which
\begin{equation}
P f \lesssim \sup _{1\leq j \leq 3 ^{d}}  P _{\mathcal D_j} f  .
\end{equation}
It therefore remains to see that the three conditions in Theorem~\ref{t:poisson} imply that
$ P _{\mathcal D}$ is bounded, for any choice of grid.

There is a Sawyer type testing theorem for dyadic positive operators \cite{09113437},
so that it suffices to verify the testing inequality
\begin{equation*}
\lVert   \mathbf 1_{\tilde Q}P _{\mathcal D} (\sigma  \mathbf 1_{Q} )\rVert _{L ^{p} (w)}
\lesssim ([w, \sigma ] _{p}  + \mathfrak P ) \sigma (Q) ^{1/p}.
\end{equation*}
as well as the dual estimate.  We have arrived at the point \eqref{e:PT}  in our proof of Theorem~\ref{t:parent}.
The remaining steps easily extend to this setting.

\mypara  One can also deduce a doubling parent testing type condition for the fractional integral operators
\begin{equation*}
T_ \alpha  f = \int f (x-y) \; \frac{dy} {\lvert  y\rvert ^{d \alpha } }, \qquad 0 < \alpha <1.
\end{equation*}

\begin{theorem}\label{t:fractional} Let $ 1<    \rho < \infty $ and $ 1< p \leq  q < \infty $.  There is a $ D>1$ so that this holds.
Let $ (w, \sigma ) $ be weights on $ \mathbb R ^{d} $.  These conditions are necessary and sufficient for $ P (\sigma   \cdot) $ to map $ L ^{p} (\mathbb R ^{d}, \sigma   )$
to $ L ^{q} (\mathbb R ^{d} , w ) $.   There is a finite constant $ \mathfrak P$ so that
\begin{enumerate}
\item   The pair of weights satisfy the $ A _{p,q}$ condition $ \sup _{Q} \frac { w (Q)^{1/q} \sigma (Q) ^{1/p'}} {\lvert  Q\rvert ^{\alpha } }
$.

\item  If $ Q\subset \mathbb R ^{d}$ is a cube for which $ \sigma (\rho Q) < D \sigma (Q)$,
then,   $ \lVert \mathbf 1_{ Q} T (\sigma \mathbf 1_{Q})\rVert _{L ^{q} (w)} \leq \mathfrak P \sigma (Q) ^{1/p}$.

\item If $ Q\subset \mathbb R ^{d}$ is a cube for which $ w (    {\rho Q}) < D w (  Q)$,
then,   $ \lVert \mathbf 1_{  Q}  T (w \mathbf 1_{  Q})\rVert _{L ^{p'} (\sigma )} \leq \mathfrak P w ( Q) ^{1/q'}$.

\end{enumerate}

\end{theorem}
There is a corresponding characterization of the weak type inequality.
The sketch of Theorem~\ref{t:poisson} applies to the Theorem above.

\mypara  It is an interesting question to gain information about the optimal choice of  $ D = D _{p, \rho }$ in Theorem~\ref{t:parent}.  We have not sought to do so, and comment briefly on the case of $ \rho =p =2$.
It is clear that $ D$ cannot be very small, because then the allowed  cubes on which one tests the norm of the maximal function are just too few, or subcritical for the pair of weights.

For instance, one knows that for weights $ w \in A_2$ that one has
\begin{equation*}
\lVert M \rVert _{L ^{2} (w) \to L ^{2} (w)} \lesssim [w] _{A_2}.
\end{equation*}
This estimate is sharp in the power of the $ A_2$ constant,
as is shown by considering the weight $ w (x) = \lvert  x\rvert ^{d - \epsilon } $ for $ 0< \epsilon < 1/2$.
One can calculate that    $ [w] _{A_2} \simeq \epsilon ^{-1} $,  with the cubes that demonstrate this being
those centered at the origin.   Note that with cube $ Q$ centered at the origin, one has
\begin{equation*}
w (Q) \simeq (\ell Q) ^{2d- \epsilon }.
\end{equation*}
It follows that the best possible choice of $ D= D _{2,2}$ would have to be of the order of $ 2 ^{2d}$.
We have shown that $ D = 2 ^{3d}$ is sufficient.

\mypara  For the values of $ 1< \rho < 2$, we are providing a very poor estimate of $ D _{\rho ,p}$.
Indeed, one would suspect that $ D _{\rho ,p} \to 1$ as $ \rho \downarrow 1$.
To show this, one would seem to need an improved notion of a shifted grids. The appendix of
\cite{MR1961195} gives one suggestion.
Similar sorts of questions have been addressed in \cites{MR3472823,MR3391907}.

\mypara  It is of  interest to extend the results of this paper to non-positive operators.
One easy remark is this.  Let  $ \mathcal D $ be a dyadic grid in $ \mathbb R ^{d}$,
and let $ \{\Delta _{Q} \;:\; Q\in \mathcal D\}$ be the associated martingale differences.
Define a martingale transform by
\begin{equation*}
T  f = \sum_{Q\in \mathcal D} \epsilon _Q \Delta _Q f, \qquad  \epsilon _Q \in \{\pm 1\}.
\end{equation*}
For the two weight inequality, one has the result of Nazarov-Treil-Volberg \cite{MR2407233} in the $ L ^2 $ case. This leads to the following sufficient conditions for a two weight inequality.

\begin{theorem}\label{t:MT} Let $ (w, \sigma )$ be weights on $ \mathbb R ^{d}$
 which satisfy the $ A_2$ condition \eqref{e:Ap}.  There is a constant $ D >1 $ so that these two conditions are   sufficient conditions for $ T (\sigma \cdot )$
 to map $ L ^{2} (\sigma )$ to $ L ^{2} (w)$:  For some finite constant $ \mathfrak P$,
\begin{enumerate}

\item For all cubes $ Q\in \mathcal D$ with $ \sigma (Q ^{(1)}) < D \sigma (Q)$,
there holds $ \lVert \mathbf 1_{Q} T (\sigma \mathbf 1_{Q})\rVert _{L ^{2} (w)} \leq \mathfrak P \sigma (Q) ^{1/2}$.

\item  The same condition above holds, with the roles of $ \sigma $ and $ w$ reversed.
\end{enumerate}
 Above, $ Q ^{(1)}$ is the dyadic parent of $ Q$.

\end{theorem}

We state this in the case of $ p=2$, as the $ L ^{p}$-case is much more complicated, see
Vuorinen \cite{MR3547711}.
The weak-type inequality for maximal truncations of martingale transforms does admit
an testing characterization. See \cite{MR2993026}*{Thm 4.3}.
One can consult \cite{MR2957550} for information about the continuous case.

\mypara Certain kinds of $ g$-functions have a two weight characterization \cite{MR3272028}.
That theorem can probably be relaxed to the current setting.   More involved would be the
weak-type estimate for maximal truncations of singular integrals, characterized in
\cite{MR2993026}*{}.

\mypara
Potentially more interesting is relaxing the testing conditions in the two weight inequality for the Hilbert transform \cite{MR3285858,MR3285857}.  
It seems very likely that such a result is true.

\bibliographystyle{alpha,amsplain}	


\end{document}